\input amstex1.tex
\documentstyle{amsppt1}
\magnification=\magstep1
\topmatter
\title Periods implying almost all periods, trees with snowflakes\\
and zero entropy maps\endtitle
\author A.~M.~Blokh \endauthor
\date{ August 1990, revised August 1991 }
\affil Institute for Mathematical Sciences, SUNY at Stony Brook \endaffil
\address{Institute for Mathematical Sciences, SUNY at Stony Brook, Stony Brook, NY 11794, USA}
\thanks { This paper was partly written while I was visiting MPI f\"ur Mathematik in Bonn and SUNY at Stony Brook; I would like to 
thank both institutions for their kind hospitality }
\abstract{Let $X$ be a compact tree and $f:X@>>>X$ be a continuous map; denote by $End(X)$ the number of endpoints of $X$ and by
$Edg(X)$ the number of edges of $X$. In Section 2 
we prove the following statements:
\roster
\item if $n>1$ is an integer with no prime divisors less than $End(X)+1$ and $f$ has a cycle
of period $n$, then $f$ has cycles of all periods greater than $2End(X)(n-1)$, and if $h(f)$ is its topological
entropy then $h(f)\ge \dfrac {\ln 2}{nEnd(X)-1}$;
\item if $0<n<End(X)+1$ and $E$ is the set of all periods of cycles of some interval map then there exists a continuous map 
$g:X@>>>X$ such that the set of all periods of cycles of $g$ is $nE\cup \{1\}$ where $nE=\{nk:k \in E\}$.
\endroster
 
	This implies that if $p$ is the least 
prime number greater than $End(X)$ and $f$ has cycles of all periods from $1$ to $2End(X)(p-1)$, then
$f$ has cycles of all periods (for tree maps this verifies the conjecture of M.~Misiurewicz, 
made in Bratislava in 1990). Combining the spectral decomposition theorem for graph maps (see [B1-B3]) with our results, we
prove the equivalence of the following statements for tree maps:
\roster
\item there exists $n$ such that $f$ has a cycle of period $mn$ for any $m$;
\item $h(f)>0$.
\endroster

	Note that the Misiurewicz conjecture and the last result are true for graph maps ([B4,B5]).

	In Section 3 we study properties of tree maps with zero entropy. Namely let $C$ be a periodic orbit
of a tree map $g:X@>>>X$; we call it {\it a snowflake}  
if it has certain properties related to those of Block's {\it simple periodic orbits} for interval maps [Bl]. 
We prove that the following statements are equivalent:
\roster
\item $h(f)=0$;
\item $(f,orb\,x)$ is a snowflake for every $x\in Per\,f$;
\item the period of every cycle of $f$ is of the form $2^lm$ where $m\le Edg(X)$ is an odd integer 
and all prime divisors of $m$ are less than $End(X)+1$.
\endroster}
\endtopmatter \baselineskip=18pt plus3pt
\document    
\heading{\bf 0. Introduction}\endheading

	Let us call one-dimensional compact branched manifolds {\it graphs}; we call them {\it trees} 
if they are connected and do not contain sets
homeomorphic to the circle. {\it In what follows we consider only continuous tree maps}. 
One of the
well-known and impressive results about dynamical properties of one-dimensional maps is the Sharkovskii theorem [S1] 
about the co-existence of periods of cycles for maps
of the real line. To formulate it let us introduce the following
{\it Sharkovskii ordering} for positive integers:
$$
3\prec5\prec7\prec\dots\prec2\cdot3\prec2\cdot5\prec2\cdot7\prec\dots\prec8\prec4\prec2\prec1\tag{$*$}
$$
Denote by $S(k)$ the set of all integers $m$ such that $k\prec m$ or
$k=m$ and by $S(2^\infty)$ the set $\{1,2,4,8,\dots\}$. Also denote by $P(\varphi)$ the set of periods of cycles of a map $\varphi$. 
\proclaim{\bf Theorem[S1]}
Let $g:\Bbb R@>>>\Bbb R$ be a continuous map. Then either $P(g)=\emptyset$ or there exists $k\in\Bbb N\cup 2^\infty$ such 
that $P(g)=S(k)$. Moreover, for any such $k$
there exists a map $g:[0,1]@>>>[0,1]$ with $P(g)=S(k)$ and there exists
a map $g_0:\Bbb R@>>>\Bbb R$ with $P(g_0)=\emptyset$.
\endproclaim
Clearly, one can apply the Sharkovskii theorem to an interval map; in this case $P(g)\neq\emptyset$ and hence  
there exists $k\in\Bbb N\cup 2^\infty$ such  
that $P(g)=S(k)$. So one can consider the Sharkovskii theorem as the first result describing possible sets of periods 
of cycles for tree maps.
Other information about these sets for tree maps is contained in [ALM] for maps of the letter Y 
and [Ba] for maps of the $n$-od.
 
	The Sharkovskii theorem implies that if a map $f:\Bbb R@>>>\Bbb R$ has a cycle of period 3 then it has cycles of all periods. 
The following conjecture, formulated by M.~Misiurewicz at the Problem Session at
Czecho-Slovak Summer Mathematical School near Bratislava in 1990, is related to the aforementioned property
of maps of the real line.
\proclaim{\bf Misiurewicz Conjecture}
For a graph $X$ there exists an integer $L=L(X)$ such that for a map $f:X@>>>X$ the
inclusion $P(f)\supset\{1,2,\dots,L\}$ implies that $P(f)=\Bbb N$. \endproclaim

	In Section 2 we verify the Misiurewicz conjecture for tree maps and get more information about 
sets of periods of cycles for tree maps. The general verification of this conjecture for arbitrary continuous
graph maps may be found in [B4,B5]. 
Note that all results of the paper are true in the same formulations for finite unions of connected trees;
the corresponding extension of our results is left to the reader.

	Fix a tree $X$. We use the terms ``vertex'', ``edge'' and ``endpoint'' in the usual sense. Denote the number of 
edges of $X$ by $Edg(X)$ and the number of endpoints of $X$ by $End(X)$. 
An integer $n$ is said to have {\it the ap-property for X} or to be {\it an ap-number for X} if there exists an integer $j(n,X)$
such that any map $f:X@>>>X$ with a cycle of period $n$ also has cycles of all periods greater than
$j(n,X)$. Theorem 1 describes the set of all ap-integers for $X$.   
\proclaim{\bf Theorem 1} Let $X$ be a tree. Then the following holds.
\roster
\item Let $n>1$ be an integer with no prime divisors less than $End(X)+1$. If a map $f:X@>>>X$ has a cycle
of period $n$, then $f$ has cycles of all periods greater than $2End(X)(n-1)$. Moreover, $h(f)\ge \dfrac {\ln 2}{nEnd(X)-1}$.
\item Let $1\le n\le End(X)$ and $E$ be the set of all periods of cycles of some interval map. Then there exists a map 
$f:X@>>>X$ such that the set of all periods of cycles of $f$ is 
$\{1\}\cup nE$, where $nE\equiv \{nk:k \in E\}$. 
\endroster \endproclaim

	Thus integers with no prime divisors less than $End(X)+1$ form the set of ap-integers for $X$. For interval maps
this implies that the set of all ap-integers coincides with the set of all odd numbers greater than $1$; clearly, this also may be 
deduced from Sharkovskii's theorem. 

	Let us formulate some corollaries of Theorem 1.
\proclaim{\bf Corollary 1 (cf. [BF])}
Let $f:X@>>>X$ be a cycle of period $n=pk$ where $p>1$ has no prime divisors less than $End(X)+1$. Then
$h(f)\ge \dfrac {\ln 2}{k[pEnd(X)-1]}>\dfrac {\ln 2}{nEnd(X)-n}$.\endproclaim

	The next corollary verifies the Misiurewicz conjecture for tree maps. 
\proclaim{\bf Corollary 2}
Let $p$ be the least prime number greater than $End(X)$. If $f:X@>>>X$ has cycles of all periods from $1$ to $2End(X)(p-1)$ then
$f$ has cycles of all periods.\endproclaim

	Theorem 1 and the spectral decomposition theorem for graph maps ([B1-B3]) imply
\proclaim{\bf Corollary 3}
The following two statements are equivalent:
\roster
\item there exists $n$ such that $f$ a cycle of period $mn$ for any $m$;
\item $h(f)>0$.
\endroster \endproclaim

	Note that in fact Corollary 3 is true for arbitrary graph maps ([B4,B5]; the different proof may be found in [LM]).

	In Section 3 we study properties of tree maps with zero entropy.
It was proved in the papers [BF], [MS] and [M] that for an interval map $g$ the fact that $h(g)=0$ is equivalent to 
$P(g)\subset \{1,2,2^2,\dots\}$. The structure of cycles for zero entropy interval maps
was studied in [Bl]. In Section 3 we generalize the above mentioned results to the case of tree maps. First we need some
definitions.

	Suppose that $A\subset X$; denote by $[A]$ 
the smallest connected set containing $A$ and 
call $[A]$ {\it the connected hull} of $A$. The definition makes sense because $X$ is a tree.
	 Let us give a definition of {\it a cycle of sets}.  Namely 
suppose that there are connected sets $\{C_i\}^{n-1}_{i=0}$ such that $fC_i\subset C_{i+1}$ for $0\le i\le n-1, fC_{n-1}\subset C_0$ and  
the sets $C_i,\;0\le i \le n-1$ are pairwise disjoint. Then we call the collection of sets 
$\{C_i\}^{n-1}_{i=0}$ {\it a cycle of sets (of period $n$)}. Note that by the definition speaking of
a cycle of sets we always mean {\it connected sets}. It is convenient to order sets in a cycle
in accordance with the way they are permuted, so usually we will index them
in this order. Finally suppose that
$Z\subset X$ and $[Z]\setminus Z$ is a connected set; then we call the set $Z$ {\it surrounding}.

	Let us pass to the definition of a snowflake. Suppose that $f:X@>>>X$ is continuous and $C\subset X$ is a connected 
invariant set. We say that {\it $f|C$ is a snowflake map (of type \linebreak $m_0=1<m_1<\dots<m_k$)} 
if there exists a collection of nested cycles of  
sets 
$$C=Y^0_0\supset \bigcup^{m_1-1}_{r=0}Y^1_r \supset \dots \supset \bigcup^{m_k-1}_{r=0}Y^k_r,\quad k>0$$ 
such that the following properties hold:

	1) $Y^0_0\supset Y^1_0\supset \dots \supset Y^k_0$, so by the definition of cycle of sets it is clear that
if $s\equiv t\pmod{m_i}$ with $0\le t<m_i,\;0\le s<m_{i+1},$ then $Y^i_t \supset Y^{i+1}_s$ (in what follows we will 
say that the sets $Y^i_r$ and the cycle of sets $\bigcup^{m_i-1}_{r=0}Y^i_r$ are {\it of level $i$});

	2) the set $\bigcup_{s\equiv t\pmod{m_i}}Y^{i+1}_s$ is surrounding for $0\le t<m_i$
(in other words, all sets of the next level belonging to some set of the previous level form a surrounding set);
in particular, $\bigcup^{m_1-1}_{r=0}Y^1_r$ is a surrounding set.

	The cycles of sets $\{\bigcup^{m_i-1}_{r=0}Y^i_r\}^k_{i=0}$ are called {\it generating}.   
Now suppose that $A$ is a cycle of sets of period $n$. Clearly, $\bigcup_{i\ge 0}f^i[A]=C$ is a connected invariant set; if  
$f|C$ is a snowflake map (of type $m_0=1<m_1<\dots<m_k=n$) and each generating set on the last level contains 
exactly one set from $A$, then we say that $(f,A)$
is {\it a snowflake (of type $m_0=1<m_1<\dots<m_k$)}. Note that in fact $C$ is the minimal by inclusion 
connected invariant set containing $A$;
the connected hull $[A]$ itself may not be invariant, so in order to get such a set we need to consider the whole orbit of $[A]$.

	Let us give some examples of snowflakes for tree maps. To begin with let us consider maps of the interval.
Then the only possible non-connected surrounding sets $Z$ are those having exactly two connected components, because
otherwise $[Z]\setminus Z$ is not connected which contradicts to the definition of surrounding sets. Now let
$f:[0,1]@>>>[0,1]$ be a continuous interval map, let $A$ be a cycle of sets and $(f,A)$ 
be a snowflake of type $m_0=1<m_1<\dots<m_k$. Set $C=\bigcup_{i\ge 0}f^i[A]$
and suppose that  $\{\bigcup^{m_i-1}_{r=0}Y^i_r\}^k_{i=0}$ are generating sets (so every $Y^i_r$ is an interval). Then by the definition 
$\bigcup^{m_1-1}_{r=0}Y^1_r$ is a surrounding set so $m_1=2$ and $Y^1_0,Y^1_1$ are simply intervals which are interchanged
by $f$. Similarly considering other levels of the snowflake $(f,A)$ we can see that
the picture on each level is close to that on the first one; in other words, for any $0\le i<k$ we have $m_{i+1}=2m_i$, 
for any $0\le t<m_i$ the intervals $Y^{i+1}_t,Y^{i+1}_{t+m_i}$ are the only intervals of level $i+1$ lying inside $Y^i_t$
and they exchange their places under the appropriate iterations of $f$ 
(namely $f^{m_i}Y^{i+1}_t\subset Y^{i+1}_{t+m_i},\;f^{m_i}Y^{i+1}_{t+m_i}\subset Y^{i+1}_t$). So we see that
the definition of a snowflake generalizes the definition of a simple periodic orbit (see [Bl]) to tree maps. 

	A similar situation takes place in the general case. Inside a set of any level the sets from the next one form
a surrounding set, so the pictures on each level inside each set are analogous; one can consider this as a sort of
self-similarity. Let us also point out that if $Z$ is a surrounding set then by the definition $A=[Z]\setminus Z$ is connected,
and geometrically $[Z]$ may be obtained by ``sticking'' components of $Z$ to the endpoints of $A$.       
 
	The most important case which we consider 
is the one when a cycle of sets $A$ is in fact an orbit of a periodic point. 
The main result is the following
\proclaim{\bf Theorem 2} Let $X$ be a tree. Then the following statements are equivalent:
\roster
\item $h(f)=0$;
\item $(f,orb\,x)$ is a snowflake for every $x\in Per\,f$;
\item every $n\in P(f)$ is of form $n=2^lm$ where $m\le Edg(X)$ is an odd integer 
and all prime divisors of $m$ are less than $End(X)+1$.
\endroster \endproclaim
\heading{\bf Notation}\endheading

	$f^n$ is the n-fold iterate of a map $f$;

	$\overline Z$ is the closure of Z;

	$\text{orb}\,x\equiv\{f^nx\}_{n=0}^\infty$ is the {\it orbit (trajectory)} of $x$;

	$\text{Per}\,f$ is the set of all periodic points of a map $f$;

	$P(f)$ is the set of all periods of periodic points of a map $f$;

	$h(f)$ is the topological entropy of a map $f$.
\heading{\bf 1. Preliminary lemmas}\endheading

	We first give some definitions. Let $X$ be a tree. By {\it an interval} we mean a homeomorphic image $h[0,1]$ of an interval
$[0,1]$ in $X$ regardless of whether it contains vertices of $X$ or not; we also consider degenerate intervals, i.e. points. 
Note that the notation we use for an interval corresponds to that from the definition of a connected hull of a set. 
Points $h(0)=a,h(1)=b$ are called  {\it endpoints} of the interval
$h[0,1] \equiv [a,b]$; clearly, there exists
a unique interval $[a,b]$ with given endpoints $a,b$. Moreover, let us denote
intervals of different types in the following way: $(a,b]\equiv[a,b]\setminus \{a\},[a,b)\equiv[a,b]\setminus \{b\},
(a,b)\equiv[a,b]\setminus \{a,b\}$. Furthermore, let $h:[0,1]@>>>[a,b],\; h(0)=a,\; h(1)=b$ be a homeomorphism; if $x\in [a,b],
y\in [a,b]$ and $h^{-1}(x)<h^{-1}(y)$ then say that {\it $x$ is closer to $a$ than $y$ 
(or $y$ is further from $a$ than $x$) on the interval $[a,b]$} (in fact
we will not mention the interval if it is clear which one we mean). Clearly, the definition is correct. 
In the similar way we will speak about subsets
of intervals in $X$; in this case by {\it $C$ is closer to $a$ than $D$ } we mean that for any $c\in C, d\in D$ either
$h^{-1}(c)<h^{-1}(d)$ or $c=d$. In what follows we consider a continuous map $f:X@>>>X$. 
	\proclaim{\bf Lemma 1}
Let $[a,b],[c,d]$ be intervals and $f[a,b]\supset[c,d],\; (fa,c)\cap(c,d)=\emptyset,\;(d,fb)\cap(c,d)=\emptyset$.
Suppose also that $I_0,I_1,\dots,I_k \subset [c,d]$ are intervals with pairwise disjoint interiors containing no
vertices of $X$ and that $I_{i+1}$ is further from $c$ than $I_i$ for $\;0\le i\le k-1$. Then there exist intervals
$J_0,J_1,\dots,J_k\subset[a,b]$ with pairwise disjoint interiors such that $J_{i+1}$ is further from $a$ than $J_i$
for $0\le i\le k-1$ and $fJ_i=I_i,\;0\le i\le k$.\endproclaim
	\demo{Proof} Clearly, for any $0\le i\le k$ there exist intervals $L\subset[a,b]$ such that $fL=I_i$. 
Indeed, let $I_i=[x,y]$ where $x$ is closer to $c$
than $y$. Choose the closest to $a$ preimage of $y$ and denote it by $y_{-1}$. Then choose the preimage of $x$ 
closest to $y_{-1}$ in $[a,y_{-1}]$,  
and denote it by $x_{-1}$. It is easy to see that $f[x_{-1},y_{-1}]=[x,y]$. Say that an interval $L$ is {\it good}
if $fL=I_i$ for some $i$ and inclusion $M\subsetneqq L$ implies that $fM\neq I_i$. Choose for $0\le i\le k$ the closest to $a$
good interval $J_i$ such that $fJ_i=I_i$. The relations $(fa,c)\cap (c,d)=\emptyset, (d,fb)\cap (c,d)=\emptyset$ easily imply now that
$J_i$ is closer to $a$ than $J_{i+1}$ for $0\le i\le k-1$ which completes the proof. \qed
	\enddemo
	\proclaim{\bf Lemma 2}
Let $J_0=[c_0,d_0],J_1=[c_1,d_1],\dots,J_k=[c_k,d_k]$ be intervals and $0=n_0<n_1<\dots <n_{k+1}$ be integers. Suppose also that 
for any $0\le i\le k-1$, we have $f^{n_{i+1}-n_i}J_i\supset J_{i+1}, (d_{i+1},f^{n_{i+1}-n_i}d_i) \cap (c_{i+1},d_{i+1})=\emptyset$
and similarly $f^{n_{k+1}-n_k}J_k\supset J_{0}, (d_{0},f^{n_{k+1}-n_k}d_k) \cap (c_{0},d_{0})=\emptyset$. Then there exists $z\in J_0$
such that $f^{n_i}z\in J_i (0\le i\le k)$ and $f^{n_{k+1}}z=z$.\endproclaim
	\proclaim{\bf Remark}{\rm In particular, if $f^{n_{i+1}-n_i}J_i\supset J_{i+1},\quad f^{n_i}d_0=d_i$ for $0\le i\le k$ and also
$f^{n_{k+1}-n_k}J_k\supset J_{0}, (d_{0},f^{n_{k+1}-n_k}d_k) \cap (c_{0},d_{0})=\emptyset$ then 
there exists $z\in J_0$ such that $f^{n_i}z\in J_i (0\le i\le k)$ and $f^{n_{k+1}}z=z$ (note that 
$(d_{0},f^{n_{k+1}-n_k}d_k)=(d_{0}, f^{n_{k+1}}d_0)$).}\endproclaim
	\demo{Proof} We divide the proof into several steps.
	\proclaim{Step 1} There exist a number $M$ and intervals $L_0,L_1,\dots,L_M \subset J_0$ such that the following holds:
\roster
\item the interiors of intervals $L_0,L_1,\dots ,L_M$ are pairwise disjoint;
\item $f^{n_i}L_j$ is an interval $(0\le i\le k+1,0\le j\le M)$;
\item $f^{n_i}L_j\subset J_i \;(0\le i\le k, 0\le j\le M)$;
\item $f^{n_{k+1}}L_0\cup f^{n_{k+1}}L_1\cup \dots \cup f^{n_{k+1}}L_M=J_0$ and 
interiors of intervals $f^{n_{k+1}}L_j, 0\le j\le M$ are pairwise disjoint;
\item for any $0\le i\le M-1$, the interval $L_i$ is closer to $c_0$ than $L_{i+1}$ and the interval $f^{n_{k+1}}L_i$
is closer to $c_0$ than $f^{n_{k+1}}L_{i+1}$.
\endroster \endproclaim

	First choose intervals $N_0,N_1,\dots,N_m$ in such a way that their union is $J_0$ and 
their interiors are pairwise disjoint and do not contain
vertices of $X$; moreover, we may assume that $N_i$ is closer to $c_0$ than $N_{i+1}$ for $0\le i\le m-1$. Furthermore, 
choose a point $x_k\in J_k$ such that $f^{n_{k+1}-n_k}x_k=c_0$. Now by Lemma 1 we can find intervals 
$T_0,T_1,\dots,T_s\subset [x_k,d_k]$ with pairwise disjoint interiors in such a way that $f^{n_{k+1}-n_k}T_i=N_i, 0\le i\le s$, and
$T_i$ is closer to $x_k$ than $T_{i+1}, 0\le i\le s-1$. It remains to divide the intervals $T_i$ into subintervals in such a way that
these new subintervals have pairwise disjoint interiors, do not contain vertices of $X$, and still are ordered in the
sense of ``closer - further'' ordering on the interval $[x_k,d_k]$. Going on with this construction and making use of
Lemma 1, we will find intervals $L_0,L_1,\dots,L_M$ with the required  properties.
	\proclaim{Step 2}
In the situation of Step 1, there exists a point $z\in \bigcup^M_{i=0}L_i$ such that $f^{n_{k+1}}z=z$.\endproclaim

	Denote $f^{n_{k+1}}$ by $g$. We may assume that $J_0=[0,1]$ and intervals $L_0,L_1,\dots ,L_M$ and 
$gL_0,gL_1,\dots ,gL_M$ increase in the usual sense. Now the fact that 
$\bigcup^M_{i=1}gL_i=[0,1] \supset \bigcup^M_{i=1}L_i$ easily implies that $\sup g|L_M=1 \ge L_M,\; \inf g|L_0=0 \le L_0$ and
so there exists $i$ such that $gL_i\supset L_i$. Indeed, we need to find $i$ such that $\sup g|L_i\ge L_i$ and
$\inf g|L_i\le L_i$. Clearly, the fact that $\inf g|L_{j+1}>L_{j+1}$ implies that $\sup g|L_j>L_j$ (because the 
intervals $\{L_j\}$ are ordered by increasing and at the same time $\bigcup^M_{i=1}gL_i=[0,1]$). Now take
the maximal $i$ such that $\inf g|L_i\le L_i$. If $i=M$ then $gL_M\supset L_M$ and we are done.
If $i<M$ then $\inf g|L_{i+1}>L_{i+1}$ and so $\sup g|L_i>L_i$ and $gL_i\supset L_i$ which completes the proof. \qed \enddemo

	A connected closed set $Y\subset X$ is called a {\it subtree}.
	\proclaim{\bf Lemma 3}
Let $X$ be a tree, $Y\subset X$ be a subtree and $f:Y@>>>X$ be a continuous map such that if $a\in Y$ then $(a,fa]\cap Y\ne \emptyset$.
Then there exists $z\in Y$ such that $fz=z$. \endproclaim
	\demo{Proof} Let us construct a map $g:X@>>>X$ in the following way. First define a map $h:X@>>>Y$ so that if $x\in Y$
then $hx=x$ and if $x\notin Y$ then $hx=y$ where $y\in Y$ is the unique point with $(y,x]\cap Y=\emptyset$. Now consider a map
$g=f\circ h:X@>>>X$. Then there exists $z\in X$ such that $gz=z$. If $z\in Y$ then $hz=z=fz$ and we are done. Let $z\notin Y$. 
Then
$hz=y$ where $y\in Y$ and $(y,z]\cap Y=\emptyset$; at the same time $gz=f(hz)=fy=z$, so $(y,fy]\cap Y=\emptyset$ which is a contradiction.\qed \enddemo
	\proclaim{\bf Lemma 4}
Let $Y\subset X$ be a subtree, $f:X@>>>X$ be a continuous map. Then there exists a point $y\in Y$ such that for any $z\in Y$ the relation
$fz\in Y$ implies the inclusion $f[y,z]\supset [y,fz]$ and either $fy=y$ or $fy\notin Y$ and $(y,fy]\cap Y=\emptyset$.\endproclaim
	\demo{Proof} Consider the case when there is no fixed point in $Y$. Then by Lemma 3 $(y,fy]\cap Y=\emptyset$ 
for some $y\in Y$. Now the properties of trees imply the conclusion.\qed\enddemo

	In what follows we call the point $y\in Y$ existing by Lemma 4 {\it the basic point for $(f,Y)$}.
\heading{\bf 1. The description of ap-numbers for tree maps}\endheading

	First we need the following definition:
if $x\in Per f$ then   
we call points $a,b\in orb\,x$ {\it neighboring} if
$(a,b)\cap orb\,x=\emptyset$. 
	\proclaim{\bf Theorem 1} Let $X$ be a tree. Then the following holds.
\roster
\item Let $n>1$ be an integer with no prime divisors less than $End(X)+1$. If a map $f:X@>>>X$ has a cycle
of period $n$ then $f$ has cycles of all periods greater than $2End(X)(n-1)$. Moreover, $h(f)\ge \dfrac {\ln 2}{nEnd(X)-1}$.
\item Let $1\le n\le End(X)$ and $E$ be the set of all periods of cycles of some interval map. Then there exists a continuous map 
$f:X@>>>X$ such that the set of all periods of cycles of $f$ is $\{1\}\cup nE$, where $nE\equiv \{nk:k \in E\}$. 
\endroster \endproclaim
	\demo{Proof} We start with statement 1). Let $x$ be a periodic point of period $n>1$ where $n$ has no prime divisors 
less than $End(X)+1$. Let $y$ be a basic point for $(f,[orb\,x])$; then $y\in [orb\,x]\setminus orb\,x$. 
Consider the connected component $Z$ of $[orb\,x]\setminus orb\,x$ such that $y\in Z$. If $z_1,z_2,\dots,z_l$ are endpoints of $Z$
then $z_i\in orb\,x$ and $(y,z_i)\cap orb\,x=\emptyset, 1\le i\le l$. Denote by $Z_i$ the connected component of the set 
$[orb\,x]\setminus Z$ containing $z_i$ and let $Y_i=Z_i\cap orb\,x$. We divide the rest of the proof into steps, but first let us make 
the following two quite simple remarks: 1) $l\le End(X)$; 2) $n\ge 3$. We also need the following easy property 
which we formulate without proof.
	\proclaim{\bf Property A} If $\{A_1,A_2,\dots,A_n\}$ are sets and $B=\bigcup^n_{i=1}[A_i]$ is connected then
$B=[\bigcup^n_{i=1}A_i]$. 
\endproclaim
	\proclaim{Step 1}
There exist two neighboring points $a,b\in orb\,x$ such that $b\in (a,y)$ and $y\in f^{l-1}(a,b)$.\endproclaim

	Let us describe the following procedure.
Let $F_1,\dots,F_m$ be pairwise disjoint subsets of $orb\,x=\bigcup^m_{i=1}F_i$ such that $[F_1],\dots,[F_m]$ are pairwise disjoint
subtrees of $X$; denote  $\bigcup^m_{i=1}[F_i]$ by $D_0$. Now consider the set $D_1=\bigcup^m_{i=1}([fF_i]\cup [F_i])$; let
$G_1,\dots,G_u$ be the connected components of $D_1$.
Denoting $H_1=G_1\cap orb\,x,\dots,H_u=G_u\cap orb\,x$, we can easily see that $G_i=[H_i],1\le i\le u$.
Indeed, denote by $\Cal A_1$ the family of all sets of type $f^rF_i, 1\le i\le m, r=1,2$. Now consider 
the set $G_j$. Then by the definition there is a subfamily $\Cal B^j\subset \Cal A_1$ such that
$G_j=\bigcup_{E\in \Cal B^j}[E],\;H_j=\bigcup_{E\in \Cal B^j}E$ and by Property A we have $G_j=[H_j]$. 
Thus the procedure
of constructing the pairwise disjoint subtrees may go on.

	Let us show that if we start the procedure in question with $m\le End(X)$ subtrees then after at most $m-1$ steps
we get the set $[orb\,x]$ (in other words we are going to show that $D_{m-1}=[orb\,x]$). Indeed, by the
properties of the number $n$ we see that $n$ and $m$ have no common divisors. Hence in the first step of the procedure  
we see that there is at least one set, say $F_1$, such that
$fF_1$ intersects with at least two of the sets $F_1,\dots,F_m$ and so the number of connected components of $D_1$ is
less than or equal to $m-1$. Repeating this argument we get the conclusion.

	It is quite easy to give the exact formula for sets $D_i$. However we need here only to show that 
$D_{j}\subset \bigcup^m_{i=1}\bigcup^{j}_{s=0}f^s[F_i]\equiv S_j$.  
Clearly, it is true for $j=0,1$.
Suppose that it is the case for some $j$;we show that $D_{j+1}\subset\bigcup^m_{i=1}\bigcup^{j+1}_{s=0}f^s[F_i]$.
Indeed, by the construction $D_{j+1}\subset D_j\cup fD_j\subset S_j\cup fS_j=S_{j+1}$ and we are done.  
	Finally we have that $[orb\,x]=D_{m-1}\subset \bigcup^m_{i=1}\bigcup^{m-1}_{s=0}f^s[F_i]$. Now let us start our procedure with
the sets $[Y_1]=Z_1,\dots,[Y_l]=Z_l$; then after $l-1$ steps we get the set $[orb\,x]$.
In other words, $[orb\,x]\subset \bigcup^l_{i=1}\bigcup^{l-1}_{s=0}f^sZ_i$. Thus there exist $s\le l-1$ and two
neighboring points $a,b\in orb\,x$ such that $b\in (a,y)$ and $y\in f^s(a,b)$; by the properties of basic points
(see Lemma 4) this implies Step 1. 

	Choose a point $\zeta \in (a,b)$ such that $f^{l-1}\zeta=y$; let for definiteness $f^{l-1}[a,\zeta]\supset[y,z_1]$. 
	\proclaim{Step 2}
There exist integers $p,q$ and $r$ such that $f^p[y,z_1]\supset [y,z_q],\,f^r[y,z_q]\supset [y,z_q]$ where
$1\le r,\,p+r\le l\le End(X)$.\endproclaim

	Consider for any $j\le l$ an integer $s(j)$ such that $[y,fz_j]\supset [y,z_{s(j)}]$. Then Lemma 4 easily implies
Step 2 ($s^p(1)=q=s^r(q)$ is an $r$-periodic point of the map $s$).

	Denote by $D$ the set $orb_s(q)=\{q,s(q),\dots,s^{r-1}(q)\}$.
	\proclaim{Step 3}
For any $v\ge (n-1)r$ and $t\in D$ we have $f^v[y,z_t]\supset [orb\,x]$.\endproclaim

	Clearly, if $B_j=f^{rj}[y,z_t]\cap orb\,x$ then $B_j\cup f^rB_j\subset B_{j+1}\;(\forall j)$. Thus 
$\bigcup^{n-1}_{j=0}f^{rj}z_t \subset f^{(n-1)r}[y,z_t]$. But $r\le End(X)$ and hence $r$ and $n$ have no common
divisors. Therefore  $\bigcup^{n-1}_{j=0}f^{rj}z_t =orb\,x$ which proves Step 3.

	Now suppose that $End(X)=c,N\ge 2c(n-1)$ and make use of Lemma 2.
Consider the following sequence of intervals and iterates of $f$ (points $\zeta, a$ have been chosen in Step 1):

	0) $J_0=[\zeta, a], n_0=0$;

	1) $J_1=[y,z_1], n_1=l-1$;

	2) $J_2=[y,z_{s(1)}], n_2=l$;

	$\vdots$

	k) $J_k=[y,z_{s^k(1)}],n_k=N-(n-1)r$ where $k=N-(n-1)r-l+2$;

	k+1) $n_{k+1}=N$.

	It is easy to see that the inequalities $n\ge 3, N\ge 2c(n-1), r\ge 1$ and $c\ge l\ge p+r$ imply that
$k=N-(n-1)r-l+2 \ge (2c-r)(n-1)-l+2 \ge 2(l+p)-l+2 \ge l$. 
Hence $s^k(1)\in D$ and by Step 3,  
$f^{(n-1)r}[y,z_{s^k(1)}]\supset [orb\,x]\supset [\zeta, a]=J_0$. So by Lemma 2 (see also Remark after Lemma 2),  
there is a point $\alpha\in [\zeta, a]$ such that $f^{n_i}\alpha\in J_i\;(0\le i\le k), f^N\alpha=\alpha$.

	Let us prove that $N$ is a period of $\alpha$. Indeed, otherwise $\alpha$ has a period $m$ which is a divisor of $N$.
Consider all iterates of $\alpha$ of type $f^{n_i}\alpha, 1\le i\le k$. 
Clearly, $\dfrac N3 \ge \dfrac {2c(n-1)}3 \ge \dfrac {4l}3 >l-1=n_1$ since $c\ge l$ and $n\ge 3$. Furthermore, 
$n_k=N-(n-1)r\ge \dfrac N2$ because $N\ge 2c(n-1) \ge 2r(n-1)$. So 
$l-1=n_1\le \dfrac N3 < \dfrac N2 \le n_k=N-(n-1)r$. At the same time, there exists $i$ such that
$n_1 \le \dfrac N3 \le mi \le \dfrac N2 \le n_k$. Hence $f^{mi}\alpha=\alpha \in [\zeta, a]$, but on the other hand,  
$f^{mi}\alpha \in \bigcup^l_{j=1}[y,z_j] \equiv S$ where $S\cap [\zeta, a]=\emptyset$. This contradiction shows that $\alpha$ 
has a period $N$. 

	To estimate $h(f)$ it is enough to observe that $f^{l-1+p+r(n-1)}[\zeta,a]\supset [\zeta, a]\cup S$,\linebreak  
$f^{l-1+p+r(n-1)}S\supset [\zeta,a]\cup S$ and at the same time $l-1+p+r(n-1)\le nc-1$. 
By usual arguments, this implies that $h(f) \ge \dfrac {\ln 2}{nc-1}$.

	Let us pass to statement 2) of Theorem 1. Let $l\le m\le End(X)$ and $g:[0,1]@>>>[0,1]$ be a map with 
$P(g)=E$. We may assume that $g(0)=0,g(1)=1$. Let us construct a map $f:X@>>>X$ such that 
$P(f)={1}\cup mE$ where $mE\equiv \{mk:k\in E\}$. 
First fix $m$ endpoints $z_1,\dots,z_m$ of $X$. For any 
$1\le i\le m$, there exists a single edge $[z_i,y_i]$ containing $z_i$; choose a point $x_i\in (z_i,y_i)$. Then 
choose $x\in X\setminus \bigcup^m_{i=1}[z_i,y_i)$. Now construct a continuous  map $f$ with the following properties. 
\item{(1)} The map $f$ outside $\bigcup^m_{i=1}[z_i,y_i)$ is identity.
\item{(2)} For any $1\le i\le m-1$ the map $f$ is injective on $[y_i,z_i]$ and maps the edge $[y_i,z_i]$ onto the union of the intervals 
$[y_i,x]\cup [x,z_{i+1}]$ in such a way that $fy_i=y_i, fx_i=x_{i+1}, fz_i=z_{i+1}$; so $f[z_i,x_i]=[z_{i+1},x_{i+1}]$.
\item{(3)} Define $f|[z_m,y_m]$ in such a way that the following holds:
\itemitem{(a)} $f|[x_m,y_m]$ is injective, $fy_m=y_m, fx_m=x_1, f[x_m,y_m]=[x_1,y_m]$;
\itemitem{(b)} $f[x_m,z_m]=[x_1,z_1]$ (which implies that $f^m[x_1,z_1]=[x_1,z_1]$) and moreover,\linebreak 
$f^m|[x_1,z_1]$ is topologically conjugate to the map $g$.

	It is easy to see now that $P(f)=\{1\}\cup mE$ where $mE\equiv \{mk:k\in E\}$. \qed \enddemo

\proclaim{\bf Corollary 1 (cf.[BF])}
Let $f:X@>>>X$ be a cycle of period $n=pk$ where $p>1$ has no prime divisors less than $End(X)+1$. Then
$h(f)\ge \dfrac {\ln 2}{k[pEnd(X)-1]}>\dfrac {\ln 2}{nEnd(X)-n}$.\endproclaim
\demo{Proof} It is enough to consider the map $f^k$ and apply Theorem 1. \qed \enddemo
\proclaim{\bf Corollary 2}
Let $p$ be the least prime number greater than $End(X)$. If $f:X@>>>X$ has cycles of all periods from $1$ to $2End(X)(p-1)$ then
$f$ has cycles of all periods.\endproclaim
\demo{Proof} The proof is left to the reader. \qed \enddemo 

	Corollary 3 follows from Theorem 1 and the spectral decomposition theorem for graph maps (see [B1-B3]).
\proclaim{\bf Corollary 3}
Let $f:X@>>>X$ be continuous. Then the following two statements are equivalent:
\roster
\item there exists $n$ such that $f$ has a cycle of period $mn$ for any $m$;
\item $h(f)>0$.
\endroster \endproclaim
\demo{Proof} Statement 1) implies statement 2) by Corollary 1. The inverse implication follows from the spectral decomposition 
theorem for graph maps (see [B1-B3]) and some properties of maps with the specification property.

	First we need the following definition: a graph map 
$\varphi:M@>>>N$ is called {\it monotone} if for any connected subset of $N$ its $\varphi$-preimage is a connected subset of 
$M$. We also need the definition of the specification property. Namely,
let $T:X \rightarrow X$ be a map of a compact infinite metric space $(X,d)$ into itself.
A dynamical system $(X,T)$ is said to have {\it the specification property} 
or simply {\it the specification} [DGS] if for any
$\varepsilon>0$ there exists such integer $M=M(\varepsilon)$ that for any $k>1$, for any $k$ points $x_1,x_2,\ldots,x_k\in X$,
for any integers $a_1\le b_1<a_2\le b_2<\ldots <a_k\le b_k$ with $a_i-b_{i-1} \ge M,\,2\le i\le k$ and for any integer $p$ with 
$p\ge M+b_k-a_1$ there exists a point $x\in X$ with $T^px=x$ such that $d(T^nx,T^nx_i)\le \varepsilon$ for $a_i\le n \le b_i,1\le i\le k$.
Maps with the specification have a lot of nice properties. The one we need may be easily obtained by methods similar to those from
[DGS], Section 21; it states that if $\psi$ is a map with the specification then there exists $N$ such that 
$P(\psi)\supset \{i:i>N\}$. 
 
	Now by the results of [B1-B3], the fact that the map $f:X@>>>X$ has a positive topological entropy 
implies that there exist a subtree $Y\subset X$, an integer $n$, a tree $Z$, a continuous map $g:Z@>>>Z$ with the specification,   
and a monotone map $\varphi:Y@>>>Z$ such that $f^nY=Y$ and $f^n|Y \circ \varphi=\varphi \circ g$ (i.e. 
$\varphi$ monotonically semiconjugates $f^n|Y$ to $g$). 
The aforementioned property of maps with the specification 
easily implies now that there exists a number $k$ such that $g$ has cycles of period $mk$ for any $m$. On the other hand, properties of
monotone graph maps and of continuous tree maps imply that then $f^n|Y$ has cycles of the same periods    
which completes the proof. \qed \enddemo
\heading{\bf 2. Trees and snowflakes}\endheading

	The main result of Section 2 is the following
\proclaim{\bf Theorem 2} Let $X$ be a tree. Then the following statements are equivalent:

	1) $h(f)=0$;

	2) $(f,orb\,x)$ is a snowflake for every $x\in Per\,f$;

	3) every $n\in P(f)$ is of form $n=2^lm$ where $m\le Edg(X)$ is an odd integer 
and all prime divisors of $m$ are less than $End(X)+1$.\endproclaim

	The definitions and notation for objects we are going to deal with in Section 2 (such as 
{\it a connected hull, a cycle of sets,
a snowflake, } etc.) may be found in the Introduction.

	First we prove Proposition 1 which is in fact a part of Theorem 2.
\proclaim{\bf Proposition 1}
Let $X$ be a tree, $h(f)=0, a\in Per\,f$ is a periodic point of period $M$. Then $(f,orb\,a)$ is a snowflake.\endproclaim
\demo{Proof}
Set $A=\bigcup_{i\ge 0}f^i[orb\,a]$ and let $y\in [orb\,a]$ be a basic point for $(f,[orb\,a])$.	
Consider the family $\Cal R$ of all cycles of sets of periods greater 
than $1$ which contain $orb\,a$ and belong to $A$. 
Let $n$ be the smallest period of a cycle of sets from $\Cal R$. Then by Zorn lemma there exists an element of $\Cal R$
which is maximal by inclusion in $\Cal R$ and has a period $n$. Clearly, $y\notin B$ 
(otherwise, the period of $B$ is $1$). Let $Z$ be a connected component of the set $A\setminus B$  containing $y$. 
We will show that $A\setminus B=Z$. First we prove the following
\proclaim{Property 1} If $C\subset A$ is connected and strictly contains some component of $B$, then 
$\bigcup_{i\ge 0}f^iC=orb\,C=A$.\endproclaim

	Indeed, if $y\in orb\,C=\bigcup^\infty_{n=0}f^nC$, then by the definition of $A$ we have $A=orb\,C$. Suppose that $y\notin orb\,C$.
Then $orb\,C \ne A$ and so $orb\,C$ is not a connected set. This easily implies that $orb\,C$ is in fact a cycle of sets
(these are exactly the components of $orb\,C$) and at the same time $orb\,C \supsetneqq B$, 
which is a contradiction and completes the proof of Property 1.

	Consider several cases. First suppose that there are components $B_i$ and $B_j$ 
of $B$ such that $D=\overline {B_i} \cap \overline {B_j} \neq \emptyset$. By the definition of cycle of sets we have 
$B_i \cap B_j =\emptyset$, so the properties of trees imply that $D$ consists of one point $x\in Per\,f$. Clearly, $x\in A$ 
(otherwise $x\in \overline A \setminus A$, i.e. $x$ is one of the endpoints of $A$ which is impossible
because $B_i \cap B_j=\emptyset$ and at the same time $\{x\}=\overline {B_i}\cap \overline {B_j}$). 
Thus $orb\,(B \cup x)=A$ by Property 1 and so $y\in orb\,x$ is a periodic point.   
Now the properties of basic points imply that
$y$ is a fixed point. Indeed, otherwise by Lemma 4 the interval $[y,fy]$ belongs to the closures
of several components of the set $B$ which is impossible. Hence $y=x$ is a fixed point and we see that 
$\{y\}=\{x\}=Z=A\setminus B$, which completes the proof of Property 1 in this case.

	Now suppose the closures of components of $B$ are pairwise disjoint. Let us prove that if $E$ is the maximal component 
of $A\setminus Z$ containing some component $F$ of $B$, then $E=F$. Indeed, 
suppose that $E\supsetneqq F$. Clearly, we may find a point $x\in E$ such that $[x,y] \cap F=[b,c]$ and $[x,y]=[x,b]\cup[b,c]\cup[c,y]$ 
where $[x,b)\cap(c,y]=\emptyset$; if we consider the natural ordering on the interval $[x,y]$ (see definitions in the beginning
of Section 1), we see that the point $x$ lies further from $y$ than the set $[x,y]\cap F=[b,c]$ on the interval $[x,y]$. 

	We are going to construct (using Property 1) a sort of ``symbolic dynamics'' for the map $f$ which guarantees that $h(f)>0$. 
Indeed, by Property 1 $orb\,(F\cup (c,y])=A$. Thus there exists a point $u\in (c,y]$ and  an integer $L$ such that $f^Lu=x$. 
It implies (by Lemma 2) that $f^L[u,y]\supset [y,x]$. Similarly, considering the set $F \cup [x,b)$ and making use of Property 1
one can find a point $v\in [x,b)$ and an integer $K$ such that $f^K[v,b]\supset [y,x]$. So we see that
\roster
\item $f^L[u,y]\supset [y,x] \supset [y,u] \cup [b,v]$;
\item $f^K[v,b]\supset [y,x] \supset [y,u] \cup [b,v]$;
\item $[b,v] \cap [y,u] =\emptyset$.
\endroster

	As usual, this implies that $h(f)>0$ which is a contradiction; thus $F=E$. Hence components of $B$ are exactly components of the set
$A\setminus Z$, i.e. the set $A\setminus B=[B]\setminus B=Z$ is connected. The cycle of sets $B_1=B=\bigcup^{n-1}_{i=0}Y^1_i$, where 
$Y^1_i$ are components of $B_1$, is of the first level in the construction of sets generating a snowflake $(f,orb\,a)$.  
Using the terminology from Section 0, we would say that $B$ is a surrounding set. Let $m_1=n,A_1=A$.

	Now set  
$A_2=\bigcup_{i\ge 0}f^{in}[orb\,a \cap Y^1_0] \subset B$. Clearly, $A_2$ is a cycle of sets of period $n$
with components belonging to the corresponding components of $B$. Consider a family
$\Cal P$ of all cycles of sets which belong to $A_2$, contain $orb\,a$ and have periods greater than $n$; then choose 
the minimal period of sets from $\Cal P$, denote it by $m_2$ and then choose the maximal (by inclusion) cycle of sets $B_2\in \Cal P$ 
with the period $m_2$.
Repeating the arguments we used finding the set $B_1$, one can easily prove that the set-theoretic difference between a component 
$G$ of $A_2$
and all components of $B_2$ belonging to $G$ is connected as it is required in the definition of a snowflake.
In other words, if $G$ is a connected component of $A_2$ then  
all the components of $B_2$ belonging to $G$ form a surrounding set and their connected hull coincides with $G$.

	Going on with the procedure of finding appropriate sets $A_i$ and $B_i$ we see that periods $m_i$
of sets $B_i$ increase strictly monotonically but cannot exceed $M$ which is a period of $a$. Hence the procedure
we have just described is finite. At the same time by the construction the procedure stops on the step $k$ if and only if $m_k=M$.
This means that a generating cycle of sets $B_k$ is of period $M$ and consists of $m_k$ components $\{Y^k_i\}^{m_k-1}_{i=0}$; each of 
these components contains exactly one point from $orb\,a$. So by the definition we see that $(f,orb\,a)$ is a snowflake, which 
completes the proof. 
 \qed\enddemo    
\proclaim{\bf Theorem 2} Let $X$ be a tree. Then the following statements are equivalent:
\roster
\item $h(f)=0$;
\item $(f,orb\,x)$ is a snowflake for every $x\in Per\,f$;
\item every $n\in P(f)$ is of the form $n=2^lm$, where $m\le Edg(X)$ is an odd integer 
and all prime divisors of $m$ are less than $End(X)+1$.
\endroster \endproclaim 
\demo{Proof}
By Proposition 1 statement 1) implies statement 2). Let us show that statement 3) follows from statement 2). 
Suppose that statement 2) holds. Let $a\in Per\,f$ have an odd period $n$;we shall show that $n$ has the required properties 
(i.e. $n\le Edg(X)$ and all prime divisors of $n$ are less than $End(X)+1$). 

	Indeed, let $(f,orb\,a)$ be a snowflake of type
$(m_0=1,m_1,\dots,m_k=n)$. Then $n=m_0\cdot (\dfrac {m_1}{m_0}) \cdot (\dfrac {m_2}{m_1}) \dots (\dfrac {m_k}{m_{k-1}})$. 
By 
definition, $\dfrac {m_i}{m_{i-1}}$ is the number of endpoints of a connected subset of $X$. Namely if $\{Y^i_j\}^{m_i-1}_{j=0}$ are
components of the generating for $(f,orb\,a)$ cycle of sets of level $i$, then the set 
$[\bigcup_{j\equiv r \pmod{m_{i-1}}}Y^i_j]\setminus \bigcup_{j\equiv r \pmod{m_{i-1}}}Y^i_j$ is connected and has  $\dfrac {m_i}{m_{i-1}}$
endpoints. So $\dfrac {m_i}{m_{i-1}}<End(X)+1$ which implies that all prime divisors of $n$ are less than $End(X)+1$.

	Let us show that there is no edge of $X$ containing more than one point of $orb\, a$. Suppose that there exist
an edge $[x,y]$ and points $a,b\in [x,y]\cap orb\,a\quad$ such that $(a,b) \cap orb\,a=\emptyset$. Let 
$orb\,a=B_k\subset B_{k-1}\subset \dots \subset B_0=\bigcup_{i\ge 0}f^i[orb\,a]$ be generating sets for $(f,orb\,a)$. 
Denote by $Y^i_0$ the component of $B_i$ containing $a$. Choose $i$ such that $b\notin Y^i_0$ and $b\in Y^{i-1}_0$. Then
the fact that $T=[\bigcup_{j\equiv 0 \pmod{m_{i-1}}}Y^i_j]\setminus \bigcup_{j\equiv 0 \pmod{m_{i-1}}}Y^i_j$ is a connected set and that
$Y^i_0\subset Y^{i-1}_0$ implies that
$T$ is a subinterval of $(a,b)$; thus $\dfrac {m_i}{m_{i-1}}=2$ is the number of endpoints of $T$.
But $n$ is odd and at the same time $n$ is a multiple of $\dfrac {m_i}{m_{i-1}}$, which is a contradiction. So there is no edge of
$X$ containing more than one point from $orb\,a$ and thus $n\le Edg(X)$.
Now suppose that $n\in P(f)$ is of form $n=2^lm$ where $m$ is an odd integer. Then we can consider the map
$f^{2^l}$ and apply to it the property
we have just proved (clearly, $f^{2^l}$ has a periodic point of odd period $m$), which completes the verification of 
statement 3). 

	Let us prove that statement 3) implies statement 1). Indeed, if $h(f)>0$ then by Corollary 3 
(or, as in the proof of Corollary 3, by the spectral decomposition theorem for graph maps [B1-B3])  
there exists $k$ such that $P(f)\supset \{ki:i>0\}$. Clearly, this contradicts statement 3) and proves that if statement 3)
holds then $h(f)=0$. This completes the proof of Theorem 2. \qed \enddemo

	Note that if $C$ is a cycle of sets, then 
one may consider the restriction $f|C$ from the view of relative positions of sets in $C$ ignoring
the behavior of the map outside or inside $C$ (one could call this approach {\it combinatorial}). Namely, suppose that 
there are a map $g$ and
a $g$-cycle of sets $C=\bigcup^{N-1}_{i=0}A_i$ (note that the map $g$ may be defined only on some part of $X$). We say that $g|C$ is
{\it a combinatorial snowflake (of type $m_0=1<m_1<\dots<m_k=N$)} if the following properties hold:
\roster
\item let $B^i_r=\kern -1.2em \displaystyle{\bigcup_{\overset {s\equiv r\pmod{m_i}} \to {0\le s<N}}\kern -1.2em A_s},\;0\le r<m_i$, 
then the sets $[B^i_r]$ are pairwise disjoint  
(here $0\le i \le k$) and we will say that sets $B^i_r$ are {\it of level $i$};
\item the set $\displaystyle{\bigcup_{\overset {s\equiv r\pmod{m_i}}\to {0\le s<m_{i+1}}}}\kern -1.2em [B^{i+1}_s],\;0\le r<m_i$ 
is surrounding for $0\le i <k$ 
(in other words, 
the set-theoretic difference between a connected hull of a set on some level and connected hulls of sets on 
the next level belonging to it is connected).
\endroster
\proclaim {\bf Remark} {\rm Obviously, for interval maps combinatorial snowflakes are exactly simple periodic orbits introduced by Block
in [Bl].}\endproclaim 

	It is easy to see that if $(f,C)$ is a snowflake then $f|C$ is a combinatorial snowflake. The following Proposition 2 shows that
from combinatorial point of view snowflakes are exactly those cycles which zero entropy tree maps may have.
\proclaim{\bf Proposition 2}
Let $A\subset X$ is finite and $g:A@>>>A$ is a map such that $g|A$ is a combinatorial snowflake. Then there exists a continuous map
$f:X@>>>X$ such that $h(f)=0$ and $f|A=g|A$.\endproclaim
\demo{Proof} Let $N=card\,A$ and $g|A$ be a combinatorial snowflake 
of type $m_0=1<m_1<\dots<m_k=N$. By definition the following properties
hold:
\roster
\item let $B^i_r=\kern -1.2em \displaystyle{\bigcup_{\overset {s\equiv r\pmod{m_i}} \to {0\le s<N}}\kern -1.2em A_s},\;0\le r<m_i$, 
then the sets $[B^i_r]$ are pairwise disjoint  
(here $0\le i \le k$) and all the sets $B^i_r$ are {\it of level $i$};
\item the set $\displaystyle{\bigcup_{\overset {s\equiv r\pmod{m_i}}\to {0\le s<m_{i+1}}}}\kern -1.2em [B^{i+1}_s],\;0\le r<m_i$ 
is surrounding for $0\le i <k$.
\endroster

	We will construct a map $f:X@>>>X$ by induction. The map $f$ will have a finite number of $f$-cycles with periods 
$m_0=1,m_1,\dots,m_k=N$ and every point of $X$ will tend to
one of these cycles. The map is already defined on the set $A=\bigcup^{N-1}_{j=0}B^k_j$;
namely it is the map $g$. So we need only to explain how to make a step in the construction, i.e. how to extend the map $f$ from 
the set $A^i=\bigcup^{m_i-1}_{j=0}[B^i_j]$ to the set $A^{i-1}=\bigcup^{m_{i-1}-1}_{j=0}[B^{i-1}_j]$. 

	Suppose that the map $f$ is defined on the set $A^i=\bigcup^{m_i-1}_{j=0}[B^i_j]$ such that $A^i$ is a cycle of the sets 
$[B^i_j]$ and the properties required 
in the previous paragraph hold, i.e. $f|A^i$ has finite number of cycles of periods $m_i,\dots,m_k$ and every point from $A^i$
tends to one of them. Consider sets $B^{i-1}_r, 0\le r< m_{i-1}$ and define the map $f$ on their connected hulls, i.e. on the sets
$[B^{i-1}_r]= \bigcup_{\overset s\equiv r\pmod{m_{i-1}} \to {0\le s<m_i}} [B^i_s]\cup Z^{i-1}_r$ 
where $Z^{i-1}_r$ is connected by the definition. 

	The map $f$ is already defined on the sets $[B^i_s]$ which form a cycle of sets and we need only to extend
the map $f$ to the union of connected sets $Z^{i-1}_r, 0\le r<m_{i-1}$ so that the sets $[B^{i-1}_r], 0\le r<m_{i-1}$ 
form a cycle of sets $A^{i-1}$, there exists only finite number of cycles of period $m_{i-1}$ belonging to 
$\bigcup^{m_{i-1}-1}_{r=0}Z^{i-1}_r$, and every $f$-orbit from $A_{i-1}$ 
which does not enter $A^i$ tends to one of these cycles of period $m_{i-1}$.
Taking into account that the map $f$ is already defined only on the endpoints of the sets $Z^{i-1}_r, 0\le r<m_{i-1}$ and that
by the definition these endpoints are mapped into $[B^{i-1}_{r+1}]\setminus Z^{i-1}_{r+1}$,  
one can easily construct the required extension of the map $f$. For the sake of completeness we give a sketch of the construction.

	1. Let us denote by $z_r$ the endpoint of the set $B^i_r$ which is common for this set and the corresponding
set $Z^{i-1}_j$ where $j\equiv r\pmod{m_{i-1}}$. Then find a point $x_r$ such that $(z_r,x_r]$ does not contain vertices of $X$ and
$(z_r,x_r]\subset Z^{i-1}_j$. Then find a point 
$y_j\in Z^{i-1}_j\setminus \bigcup_{r\equiv j\pmod{m_{i-1}}}(z_r,x_r]=P_j$
for any $0\le j<m_{i-1}$. 

	2. Set $f(P_j)=y_{j+1} (0\le j\le m_{i-1}-2), f(P_{m_{i-1}-1})=y_0$.

	3. Define $f|[z_r,x_r], 0\le r<m_i$ in such a way that $f[z_r,x_r]=[z_{r+1},y_{r+1}], f|[z_r,x_r]$ is injective
and if $D_r\subset [z_r,x_r]$ consists of all points $\zeta$ such that $f^{m_i}(\zeta) \in (z_r,x_r]$ then $D_r$ is an
interval, $f^{m_i}D_r=(z_r,x_r]$, there is only one periodic point belonging to $(z_r,x_r]$ and this point is of period $m_i$.

	It is easy to check that the construction in question is possible and that this way we will construct
a map with the required properies. This completes the proof of Proposition 2. \qed\enddemo

	Using methods similar to those from the proof of statement 2) from Theorem 1 or Proposition 2 
one can easily prove the following
\proclaim{\bf Proposition 3}
If $m\le End(X)$ and $k\ge 0$ then there exists a continuous map $f:X@>>>X$ and an $f$-periodic point $a$ such that $h(f)=0$ and
the period of $a$ is $2^km$.\qed \endproclaim 
\bigskip
\font\fontA=cmcsc10
{\fontA Acknowledgments.}
I would like to thank E.~M.~Coven for informing me about Baldwin's paper [Ba] and J.~Milnor and S.~Sutherland 
for looking through this paper and providing useful comments.   

\Refs
\ref \key[{\bf ALM}] \by L. Alsed\'a, J. Llibre, M. Misiurewicz \paper Periodic orbits of maps of Y 
\jour Trans. Amer. Math. Soc. \vol 313 \yr 1989 \pages 475--538 \endref
\ref \key[{\bf Ba}] \by S. Baldwin \paper An Extension of Sharkovskii Theorem to the n-od
\jour Erg. Th. and Dyn. Syst., \#5, \vol 11 \yr 1991 \pages 249-271\endref
\ref \key[{\bf Bl}] \by L. Block \paper Simple Periodic Orbits of Mappings of the Interval
\jour Trans. Amer. Math. Soc. \vol 254 \yr 1979\pages 391--398\endref
\ref \key[{\bf B1}] \by A. M. Blokh \paper On the limit behavior of one-dimensional dynamical systems.1,2 {\rm (in Russian) }\linebreak
\jour Preprints VINITI \#1156-82, \#2704-82, Moscow \yr 1982\endref
\ref \key[{\bf B2}] \by A. M. Blokh \paper Decomposition of Dynamical Systems on an Interval
\jour Russ.Math.Surv., \#5, \vol 38 \yr 1983 \pages 133--134\endref
\ref \key[{\bf B3}] \by A. M. Blokh \paper On Dynamical Systems on One-Dimensional Branched Manifolds.1 {\rm( in Russian )}
\jour Theory of Functions, Functional Analysis and Applications, Kharkov, \vol 46 \yr 1986 \pages 8--18
\moreref \paper 2 \jour Theory of Functions, Functional Analysis and Applications, Kharkov,\vol 47\yr 1986\pages 67--77
\moreref \paper 3 \jour Theory of Functions, Functional Analysis and Applications, Kharkov,\vol 48\yr 1987\pages 32--46\endref
\ref \key [{\bf B4}] \by A. M. Blokh \paper The Spectral Decomposition, Periods of Cycles and Misiurewicz Conjecture for Graph Maps
\jour (1990, submitted to ``Proceedings of the Conference on Dynamical Systems in G\"ustrow'', to appear in
Lecture Notes in Math.) \endref
\ref \key [{\bf B5}] \by A. M. Blokh \paper On Some Properties of Graph Maps: Spectral Decomposition, Misiurewicz Conjecture and Abstract
Sets of Periods \jour Max-Planck-Institut f\"ur Mathematik, Preprint \#35, June \yr 1991 \endref
\ref \key [{\bf BF}] \by R. Bowen, J. Franks \paper The Periodic Points of Maps of the Disk and the Interval 
\jour Topology \vol 15 \yr 1976 \pages 337--342 \endref 
\ref \key[{\bf DGS}] \by M. Denker, C. Grillenberger, K. Sigmund \book Ergodic Theory on Compact Spaces
\publ Lect.Notes in Math., vol.527 , Springer \publaddr Berlin \yr 1976 \endref
\ref \key[{\bf LM}] \by J. Llibre, M. Misiurewicz \paper Horseshoes, Entropy and Periods for Graph Maps \jour to appear\endref
\ref \key[{\bf S1}] \by A.N.Sharkovskii \paper Coexistence of Cycles of a Continuous Map of a Line into itself
\jour Ukr. Math. \linebreak Journal \vol 16 \yr 1964\pages 61--71 \endref
\enddocument